\documentclass{article}
\usepackage{amsfonts, amsthm, amsmath, amssymb}

\PassOptionsToPackage{hyphens}{url}\usepackage[colorlinks,citecolor=blue,bookmarks=true]{hyperref}

\theoremstyle{theorem}
\newtheorem{theorem}{Theorem}

\theoremstyle{definition}

\newcommand{\B}[1]{\mathbf{{#1}}}

\newcommand{\F}{\mathbb{F}}
\newcommand{\FF}{\mathbb{F}}

\newcommand{\NN}{\mathbb{N}}

\newcommand{\RR}{\mathbb{R}}

\newcommand{\K}{\mathbb{K}}
\newcommand{\C}{\mathcal{C}}

\newcommand{\sub}{\subseteq}
\newcommand{\sm}{\setminus}

\DeclareMathOperator{\I}{I}
\DeclareMathOperator{\Z}{Z}

\renewcommand{\L}{\mathbb{L}}

\date{}
\begin{document}

\title{Sharp Effective Finite-Field Nullstellensatz}

\author{Guy Moshkovitz and Jeffery Yu}




\maketitle

\begin{abstract}
	The (weak) Nullstellensatz over finite fields says that if $P_1,\ldots,P_m$ are $n$-variate degree-$d$ polynomials with no common zero over a finite field $\F$ then there are polynomials $R_1,\ldots,R_m$ such that
	$R_1P_1+\cdots+R_mP_m \equiv 1$.
	Green and Tao (\cite{GreenTao09}, Proposition~9.1)
	used a regularity lemma to obtain an effective proof, showing that the degrees of the polynomials $R_i$
	can be bounded independently of $n$,
	though with an Ackermann-type dependence on the other parameters $m$, $d$, and $|\F|$.
	In this paper we use the polynomial method to give a proof with a degree bound of $md(|\F|-1)$.
	We also show that the dependence on each of the parameters is the best possible up to an absolute constant.
	We further include a generalization, offered by Pete L.~Clark,
	from finite fields to arbitrary subsets in arbitrary fields, provided the polynomials $P_i$ take finitely many values on said subset.
\end{abstract}

Hilbert's Nullstellensatz is a linchpin of algebraic geometry, providing a bridge between algebra (ideals in polynomial rings) and geometry (solution sets of polynomial equations) when the underlying field is algebraically closed.
Koll\'{a}r~\cite{Kollar88} famously obtained an effective proof of Hilbert's Nullstellensatz with a sharp bound on the degrees of the polynomials involved.
When the underlying field is finite, Green and Tao~\cite{GreenTao09} proved an effective analogue of Hilbert's Nullstellensatz, but the bound furnished by their proof is extremely poor.
In this paper we obtain new, much better bounds for the finite-field analogue of Hilbert's Nullstellensatz.

\section{Upper bound.}


Two polynomials $P,Q \in \F[\B{x}]:=\F[x_1,\ldots,x_n]$ over a finite field $\F$ are functionally equal, denoted $P \equiv Q$, if $P(\B{x}) = Q(\B{x})$ for every $\B{x} \in \F^n$; put differently, $P \equiv Q$ if and only if $P = Q + C$ for some polynomial $C \in \F[\B{x}]$ in the ideal $\langle x_1^{|\F|}-x_1, \ldots, x_n^{|\F|}-x_n \rangle$.
%
%
We denote by $\F_{\le d}[\B{x}] = \{P \in \F[\B{x}] \mid \deg P \le d\}$ 
the set of polynomials over a field $\F$ of (total) degree at most $d$,
and we denote by $\Z_\F(P_1,\ldots,P_m)= \{ \B{x} \in \FF^n \mid P_1(\B{x})=\cdots=P_m(\B{x})=0\}$ the set of zeros over $\F$ (not over the algebraic closure!) of the polynomials $P_1,\ldots,P_m \in \F[\B{x}]$.


Over an algebraically closed field $\K$, Hilbert's Nullstellensatz~\cite{Hilbert1893} says that
if polynomials $P_1,\ldots,P_m,Q \in \K_{\le d}[\B{x}]$ satisfy $\Z_\K(P_1,\ldots,P_m) \sub \Z_\K(Q)$, then there exist polynomials $R_1,\ldots,R_m \in \K_{\le D}[\B{x}]$ such that $Q^{D'} = \sum_{i=1}^m R_iP_i$ for some finite $D,D' \in \NN$.\footnote{In terms of ideals, $\I(\Z_\K(P_1,\ldots,P_m)) = \sqrt{\langle P_1,\ldots,P_m \rangle}$.}
A result of Hermann~\cite{Hermann26} gives an effective proof of Hilbert's Nullstellensatz with a double-exponential bound on $D,D'$.
This bound was greatly improved by Brownawell~\cite{Brownawell87}, and an influential work of Koll\'{a}r~\cite{Kollar88} 
finally achieved a sharp bound, showing that $D,D'$ can be bounded by roughly $d^m$, and in particular, be bounded independently of $n$.


Over a finite field $\FF$, for an analogue of Hilbert's Nullstellensatz to hold, one has to replace polynomial equality with functional equality.
It is easy to see that $\Z_\F(P_1,\ldots,P_m) \sub \Z_\F(Q)$ implies the existence of \emph{functions} $R_i\colon \F^n\to\F$ satisfying $Q \equiv \sum_{i=1}^m R_iP_i$ (no need for a radical).
Since every function over a finite field can be written as a polynomial, this trivially implies the degree bound $\deg(R_i) \le n(|\F|-1)$ (recall $x_i^{|\F|}\equiv x_i$, so the degree of each variable can be assumed to be at most $|\F|-1$).
Green and Tao, in a seminal work on the relation between structure and randomness of polynomials, 
devised a certain regularity lemma for polynomials over finite fields
and applied it to obtain (among other applications) an effective version of the Nullstellensatz over finite fields~(\cite{GreenTao09}, Proposition~9.1).
Specifically, they showed that if $\F$ is a prime finite field
and $P_1,\ldots,P_m,Q \in \F_{\le d}[\B{x}]$ are polynomials of degree at most $d < |\F|$ satisfying $\Z_\F(P_1,\ldots,P_m) \sub \Z_\F(Q)$, then 
$Q \equiv \sum_{i=1}^m R_iP_i$ holds for some polynomials $R_1,\ldots,R_m \in \F_{\le D}[\B{x}]$ of degree at most some $D=D(m,d,|\F|)$ that does not depend on $n$. Unfortunately, the dependence of $D$ on the other parameters, $m$, $d$, and $|\F|$, is atrocious, being some Ackermann-type function.\footnote{See, e.g., the remark following Lemma 2.4 in~\cite{GreenTao09} for a discussion on the Ackermann-type dependence. Let us also mention that the dependence on the field size in Proposition~9.1 of~\cite{GreenTao09} is implicit, as the finite field is fixed at the beginning of that paper (e.g. $\F=\FF_2$ or $\F=\F_3$). That a dependence on the field size is necessary is shown in Section~\ref{sec:LB} below.}

We use the polynomial method to give a sharp effective proof of the finite-field Nullstellensatz (which is also short and elementary).

\begin{theorem}[Sharp finite-field Nullstellensatz]\label{theo:main}
	Let $\F$ be a finite field.
	If polynomials $P_1,\ldots,P_m,Q \in \F_{\le d}[\B{x}]$ 
	satisfy $\Z_\F(P_1,\ldots,P_m) \sub \Z_\F(Q)$, 
	then $Q \equiv \sum_{i=1}^m R_iP_i$ for some polynomials $R_1,\ldots,R_m \in \F_{\le D}[\B{x}]$ with $D = md(|\F|-1)$.
%
\end{theorem}
%
We note that an immediate corollary 
 is a finite-field analogue of the \emph{weak} Nullstellensatz, obtained by taking $Q=1$ so that the assumption is $\Z_\F(P_1,\ldots,P_m) = \emptyset$.\footnote{For $Q=1$ the proof gives the slightly better bound $md(|\F|-1) - d$; see~(\ref{eq:I_i}).}
We also note that the proof of Theorem~\ref{theo:main} has a generalization to arbitrary fields; see Section~\ref{sec:finite-null}.

Theorem~\ref{theo:main} is proved roughly as follows. Over a finite field $\F$, the indicator function for non-membership in the zero set $\Z_\F(P_1,\ldots,P_m)$ can be written as a polynomial $I$ in the ring $\F[P_1,\ldots,P_m]$ of degree at most $m(|\F|-1)$, which moreover has a zero constant term. 
The Nullstellensatz assumption implies that $Q \equiv QI$, 
and so we can write $Q$ as a combination of $P_1,\ldots,P_m$. Since the coefficients of this combination are polynomials in $\F[P_1,\ldots,P_m,Q]$, the desired degree bound follows.
The underlying idea is similar to some other applications of the polynomial method.  Broadly speaking, the polynomial method analyzes an object of interest by studying the algebraic properties of a polynomial associated with it.
In our case, the object of interest is a zero set over a finite field, 
which we analyze via a polynomial vanishing on it;
this particular idea can be traced back at least to Chevalley~\cite{Chevalley35} (see also Ax's proof of the Chevalley-Warning Theorem~\cite{Ax64}). A more recent relevant example of the polynomial method is the solution of the cap-set problem by Ellenberg-Gijswijt~\cite{EllenbergGi17} following Croot-Lev-Pach~\cite{CrootLePa17}.
The main idea in our proof is that non-membership in the zero set of $m$ polynomials over a finite field can itself be expressed as a polynomial which lies in the ideal generated by the $m$ polynomials.
The crucial point is that the number of inputs to the non-membership function is the number $m$ of polynomials rather than the number $n$ of variables of these polynomials.

%

\begin{proof}
	Put $q=|\F|$, 
	and recall that for 
	every $x \in \F$ we have $x^{q-1} \equiv [x \neq 0]$.\footnote{The Iverson bracket $[P]$ is $1$ if the predicate $P$ is true, and $0$ otherwise.}
	Let $I\colon \F^n\to\{0,1\}$ be the indicator function
	$I(\B{x}) = [\B{x} \notin \Z_\F(P_1,\ldots,P_m)]$.
	We claim that 
	$I \equiv \sum_i I_i P_i$,
	where each $I_i\in \F[\B{x}]$ is explicitly given by
	\begin{equation}\label{eq:I_i}
		I_i = P_i^{q-2}\prod_{j=i+1}^m (1-P_j^{q-1}).
	\end{equation}
	Indeed, 
	\begin{align*}
		\sum_{i=1}^m I_i P_i
		&= \sum_{i=1}^m P_i^{q-1}\prod_{j=i+1}^m (1-P_j^{q-1})
		= \sum_{i=1}^m \sum_{\substack{J \sub [m]\colon\\\min J=i}}\, -\prod_{j \in J} (-P_j^{q-1})\\
		&= -\sum_{\substack{J \sub [m]\colon\\J \neq \emptyset}}\, \prod_{j \in J} (-P_j^{q-1})
		= -\Big( \prod_{j=1}^m (1-P_j^{q-1}) - 1 \Big)
		\equiv I .
	\end{align*}
%
	Using the statement's assumption $\Z_\F(P_1,\ldots,P_m) \sub \Z_\F(Q)$ we deduce that 
	$$Q \equiv QI
	\equiv \sum_{i=1}^m (QI_i) P_i.$$
	Each $R_i:=QI_i$ can be viewed, by~(\ref{eq:I_i}), as a polynomial in the ring $\F[P_1,\ldots,P_m,Q]$, which is thus of degree at most $(m+1)(q-1)$; in fact, the explicit definition in~(\ref{eq:I_i}) implies the slightly better bound $m(q-1)$. 
	It follows that $R_i$, viewed as a polynomial in $\F[\B{x}]$, is of degree at most $md(q-1)=D$, thus completing the proof.
%
%
\end{proof}

\section{Lower Bounds.}\label{sec:LB}

We next show that a dependence in Theorem~\ref{theo:main} of the degree bound $D$ on each of the parameters $|\F|$, $d$, and $m$ is inevitable, even for the weak form of the finite-field Nullstellensatz.


For the dependence of $D$ on the field size $|\F|$, 
let us consider the polynomial $P=x^2+1 \in \F[x]$ where $\F=\F_q$ is any prime finite field with the prime number $q$ congruent to $3$ modulo $4$.
Then $\Z_\F(P) = \emptyset$ by Euler's criterion~\cite{Euler50,Euler61}, 
and so trivially $\Z_\F(P) \sub \Z_\F(1)$. However, if $1 \equiv RP$ for some $R \in \F[x]$ then, as we show next, $\deg(R) \ge q-1$.
First, note that $R \equiv P^{q-2}$ since $R(x) = 1/P(x)$ for every $x \in \F$.
Next, note that for every $i \not\equiv 0 \pmod{q-1}$ we have over $\F$ the functional equality of polynomials $x^i \equiv x^{i \pmod{q-1}}$; indeed,
for all integers $d,i \ge 1$,
$$x^{d(q-1) + i} = x^{q + (d-1)(q-1) + i-1} \equiv x \cdot x^{(d-1)(q-1) + i-1} = x^{(d-1)(q-1)+i} ,$$
as none of the exponents is negative.
Denote
$t = q-2$ and $b = \frac12(q-1) \in \NN$, so
$$R \equiv (x^2+1)^t = \sum_{i=0}^t \binom{t}{i} x^{2i}
\equiv \binom{t}{b}x^{2b} + \sum_{j=1}^{b-1} \bigg(\binom{t}{j}+\binom{t}{j+b}\bigg)x^{2j} + 1$$
as $x^{2j} \equiv x^{2(j+b)}$ for $j \not\equiv 0 \pmod{b}$.
Note that the coefficient $\binom{t}{b}$ is $\binom{q-2}{\frac12(q-1)} \not\equiv 0$ $\pmod q$, as $q$ is prime and so does not divide the numerator of the binomial coefficient.
Since every exponent of $x$ in the right hand side of the identity above is at most $2b=q-1$, and since the coefficient of $x^{q-1}$ is nonzero in $\F$, we deduce that $\deg(R) \ge q-1 = |\F|-1$, as claimed.

Note, however, that since the example above assumes that $\F$ is a prime finite field, 
it could be that $D$ depends on just 
the characteristic of $\F$ rather than its size $|\F|$. 
Nevertheless, with a little effort one can extend the example to non-prime fields: simply let $\F=\F_q$ with $q=p^k$, where $k \ge 1$ is any odd integer and $p$ is a prime congruent to $3$ modulo $4$, and take $P = x^2+1 \in \F[x]$ as before.
Then $\Z_\F(P) = \emptyset$ again; indeed, if $\sqrt{-1} \in \F_q$ then, since $\sqrt{-1} \notin \F_p$, the field extension $\L/\F_p$ with $\L=\F_p(\sqrt{-1})$ has degree $[\L : \F_p]=2$, which is impossible as $k=[\F_q : \F_p] = [\F_q : \L] \cdot [\L : \F_p]$ is odd.
Since we still have the identity $x^q \equiv x$, the only part of the argument above that remains to be checked is that the leading coefficient is nonzero in $\F$, that is, $\binom{q-2}{\frac12(q-1)} \not\equiv 0$ $\pmod p$.
Here we use a corollary of another classical result, Lucas's theorem~\cite{Lucas78} (though it would suffice to use the  earlier Kummer's theorem~\cite{Kummer52}, which itself follows from Legendre's formula):
for every prime $p$ and binomial coefficient $\binom{n}{m}$, we have $\binom{n}{m} \not\equiv 0 \pmod p$ if and only if $m_i \le n_i$ for every $i$,
where $n = \sum_{i \ge 0} n_i p^i$ and $m = \sum_{i \ge 0} m_i p^i$ are the expansions in base $p$ of $n$ and $m$, respectively.
For our binomial coefficient $\binom{q-2}{\frac12(q-1)}$ we have the following expansions in base $p$, where we recall that $q=p^k$:
$$p^k - 2 = \sum_{i=1}^{k-1} (p-1)p^i + (p-2)p^0
\quad\text{ and }\quad
\frac12(p^k-1) = \sum_{i=0}^{k-1} \Big(\frac12(p-1)\Big)p^i .$$
Since $\frac12(p-1) \le p-2$ for every $p \ge 3$, it follows from the above that $\binom{q-2}{\frac12(q-1)} \not\equiv 0 \pmod p$, 
as needed.

For the dependence of $D$ on the degree bound $d$, 
consider instead the polynomial $P=H^2+1$ where $H \in \F[x_1,\ldots,x_n]$ is any polynomial with individual degrees $\deg_{x_i}(H)=1$ (e.g., the elementary symmetric polynomial $H(\B{x})=\sum_{I \in \binom{[n]}{k}} \prod_{i \in I} x_i$ of any degree $1 \le k \le n$).
As before, $\emptyset = \Z_\F(P) \sub \Z_\F(1)$;
and if $1 \equiv RP$ then, as is not hard to check, the argument above implies that $\deg(R) \ge \deg(H)(|\F|-1) = \frac12\deg(P)(|\F|-1)$. 

Finally, the dependence of $D$ on the number $m$ of polynomials follows from results in theoretical computer science.
The Nullstellensatz over finite fields is also studied in the area of \emph{proof complexity}~\cite{BeameImKrPiPu96}, where an identity of the form
\begin{equation}\label{eq:identity}
	1=\sum_{i=1}^m R_iP_i+\sum_{i=1}^n C_i(x_i^2-x_i)
\end{equation} 
for $n$-variate polynomials $R_1,\ldots,R_m, C_1,\ldots,C_n \in \F[x_1,\ldots,x_n]$ over any finite field $\F$ is known as a \emph{Nullstellensatz refutation} for $P_1,\ldots,P_m$ (as it refutes the satisfiability of the polynomial equations $P_1(\B{x})=\cdots=P_m(\B{x})=0$ with $\B{x} \in \{0,1\}^n$).
In that context, the number $m$ of polynomials $P_i$ typically grows with the number of variables $n$, and $m$-tuples of polynomials $(P_1,\ldots,P_m)$ for which $\max_i\deg(R_i) = \Omega(n)$
 are known to exist, over every finite field~\cite{AlekhnovichRa03,BenIm10,BussGrImPi01,Grigoriev98}.
In the setting of Theorem~\ref{theo:main}, however, the crucial point is that the number of polynomials $m$ is a given parameter while the number $n$ of variables may go to infinity.
Nevertheless, as we next explain, the proof complexity lower bounds easily give polynomials in an arbitrarily large number $N$ of variables, yielding a lower bound for Theorem~\ref{theo:main} of $D \ge \Omega(m)$.
Indeed, the aforementioned lower bounds say that for any $\F$, any $m$, and \emph{some} $n=\Theta(m)$, there are polynomials $P_1,\ldots,P_{m} \in \F[x_1,\ldots,x_n]$ such that the identity~(\ref{eq:identity})
implies $\deg(R_i) = \Omega(m)$ for every $1 \le i \le m$.
Now, for any $\F$, any $m$, and any $N \ge n$, let $P_1,\ldots,P_m$ be the polynomials from before but viewed as elements of
the larger polynomial ring $\F[x_1,\ldots,x_n,x_{n+1},\ldots,x_N]$.
If $1 \equiv \sum_{i=1}^m R'_iP_i$, meaning 
$$1=\sum_{i=1}^m R'_iP_i+\sum_{i=1}^N C'_i(x_i^{|\F|}-x_i)$$ 
for some polynomials $R'_1,\ldots,R'_m,C'_1,\ldots,C'_N \in \F[x_1,\ldots,x_N]$, then by setting $x_i=0$ for every $n+1 \le i \le N$ and factoring $x_i^{|\F|}-x_i = S(x_i)(x_i^2-x_i)$ with $S(x) := \sum_{j=0}^{|\F|-2}x^j$ we get 
$$1=\sum_{i=1}^m R'_i(x_1,\ldots,x_n,0,\ldots,0)P_i+\sum_{i=1}^n C'_i(x_1,\ldots,x_n,0,\ldots,0)S(x_i)(x_i^2-x_i).$$
This identity is of the form~(\ref{eq:identity}), so $\deg(R'_i) \ge \deg(R'_i(x_1,\ldots,x_n,0,\ldots,0)) = \Omega(m)$ for every $1 \le i \le m$, as claimed.


\section{``Finite'' Nullstellensatz.}\label{sec:finite-null}

Following the initial publication of this work, Pete L.~Clark~\cite{Clark21} observed that the proof of Theorem~\ref{theo:main} extends from finite fields to arbitrary subsets of arbitrary fields, as long as the image of the polynomials on said subset is finite.
We give below a variant of this argument, included with Clark's permission.
For polynomials $P,Q$ and a set $X$,
we write $P(X) = \{P(x) \mid x \in X\}$ for the image of $P$ on $X$;
we write $P \equiv_X Q$ if $P(\B{x}) = Q(\B{x})$ for every $\B{x} \in X$, or equivalently, $P \equiv Q \pmod{\I(X)}$.\footnote{$\I(X) = \{f \in \F[x_1,\ldots,x_n] \mid \forall \B{x} \in X \colon f(\B{x}) = 0\}$ is the ideal of polynomials vanishing on $X \sub \F^n$.}


\begin{theorem}[Sharp ``Finite'' Nullstellensatz]\label{theo:finite}
	Let $\F$ be any field, $X \sub \F^n$ any set, 
	and $P_1,\ldots,P_m,Q \in \F_{\le d}[\B{x}]$ $n$-variate polynomials
	satisfying $|P_i(X)| \le F$ for every $1 \le i \le m$.
	If $\Z(P_1,\ldots,P_m) \cap X \sub \Z(Q)$
	then $Q \equiv_X \sum_{i=1}^m R_iP_i$ for some $R_1,\ldots,R_m \in \F_{\le D}[\B{x}]$ with 
	$$D = 
md\cdot\begin{cases}
	F-1	& m \ge 2\\
	F		& m=1
\end{cases}.$$
\end{theorem}

Theorem~\ref{theo:finite} in particular implies a qualitative statement proved in a paper of Clark (\cite{Clark14}, Theorem~7, ``Finitesatz''):
For any field $\F$, if $X \sub \F^n$ is finite 
then $\I(\Z(P_1,\ldots,P_m) \cap X) = \langle P_1,\ldots,P_m \rangle + \I(X)$ (without a radical, unlike Hilbert's Nullstellensatz).
In fact, Theorem~\ref{theo:finite} implies the same conclusion for infinite $X$, as long as the images $P_i(X)$ are finite; this
qualitative statement already appears to be new.

%
%
%
%

\begin{proof}
	First, suppose $\Z(P_i) \cap X = \emptyset$ for some $1 \le i \le m$.
	By interpolation, there is a univariate polynomial $I_i \in \F[x]$
	with $\deg(I_i) \le |P_i(X)|-1$ satisfying $I_i(x) = 1/x$ for every $x \in P_i(X)$ (which is well defined as $0 \notin P_i(X)$).
	Then the polynomial $R_i := Q \cdot (I_i \circ P_i) \in \F[\B{x}]$ satisfies $Q \equiv_X R_iP_i$ as needed;
	indeed, this is because $R_iP_i(\B{x}) = Q(\B{x}) \cdot I_i(P_i(\B{x})) \cdot P_i(\B{x}) = Q(\B{x})$ for every $\B{x} \in X$.
	We have 
	\begin{equation}\label{eq:bd-m1}
	\deg(R_i) \le \deg(Q) + \deg(P_i)(F-1).
	\end{equation}
	Thus, $\deg(R_i)
	\le d + d(F-1) = dF$.
	Note that $dF \le D$; indeed, 
	for $m\ge 2$ we have $dF \le 2d(F-1)$ if $F \ge 2$, whereas if $F \le 1$ then the statement is clearly true as $\deg(R_i)=0$ suffices.
	This completes the proof in this case.
	
%
%
%
%
%
%
%
%
%
%
%
%
%
	
	Henceforth, assume $\Z(P_i) \cap X \neq \emptyset$ for every $1 \le i \le m$.
	For a set $Y \sub \F$, let $\C_Y \in \F[x]$ be the univariate polynomial
	$$\C_Y(x) = 1 - \frac{\prod_{y \in Y\sm\{0\}} (y-x)}{\prod_{y \in Y\sm\{0\}} y}.$$
	Note that $\C_Y(x) = [x \neq 0]$ for every $x \in Y$,
	and $\deg(\C_Y) \le |Y \setminus \{0\}|$.
	Let $\hat{P_j} := \C_{P_i(X)} \circ P_i \in \F[\B{x}]$, 
	and note that
	\begin{equation}\label{eq:ol}
		\forall \B{x} \in X \colon \hat{P_i}(\B{x}) = [P_i(\B{x}) \neq 0]
		\quad\text{ and }\quad \deg(\hat{P_i}) \le \deg(P_i)(|P_i(X)|-1),
	\end{equation}
	where the inequality uses the assumption that $0 \in P_i(X)$.
	For every $1 \le i \le m$, let $I_i \in \F[\B{x}]$ be the polynomial
	\begin{equation}\label{eq:I_i-B}
		I_i = \frac{\hat{P_i}}{P_i} \prod_{j=i+1}^m (1-\hat{P_j}).
	\end{equation}
	Observe that $\hat{P_i}/P_i$ (and thus $I_i$) is indeed a polynomial since $\C_Y(x)/x \in \F[x]$ is, which follows from the fact that the univariate polynomial $\C_Y$ satisfies $\C_Y(0)=0$.
	Let the function $I \colon X \to \{0,1\}$ be given by $I(\B{x}) = [\B{x} \notin \Z(P_1,\ldots,P_m)]$.
	We claim that $I \equiv_X \sum_{i=1}^m I_iP_i$.
	Indeed,
	\begin{align*}
		\sum_{i=1}^m I_i P_i
		&= \sum_{i=1}^m \hat{P_i}\prod_{j=i+1}^m (1-\hat{P_j})
		= \sum_{i=1}^m \sum_{\substack{J \sub [m]\colon\\\min J=i}}\, -\prod_{j \in J} (-\hat{P_j})\\
		&= -\sum_{\substack{J \sub [m]\colon\\J \neq \emptyset}}\, \prod_{j \in J} (-\hat{P_j})
		= -\Big( \prod_{j=1}^m (1-\hat{P_j}) - 1 \Big)
		\equiv_X I,
	\end{align*}
	where the last step follows as $\prod_{j=1}^m (1-\hat{P_j}(\B{x})) = \prod_{j=1}^m [P_j(\B{x}) = 0] = [\B{x} \in \Z(P_1,\ldots,P_m)]$ for every $\B{x} \in X$.
	
	Using the statement's assumption $\Z(P_1,\ldots,P_m) \cap X \sub \Z(Q)$ we deduce that 
	$$Q \equiv_X QI
	\equiv_X \sum_{i=1}^m (QI_i) P_i.$$
	Each $R_i:=QI_i \in \F[\B{x}]$ is, by~(\ref{eq:ol}) and~(\ref{eq:I_i-B}), of degree
	$$\deg(R_i) \le \deg(Q) + \deg(P_i)(|P_i(X)|-2) + \sum_{j=i+1}^m \deg(P_j)(|P_j(X)|-1) .$$
	In particular, 
	$\deg(R_i) \le md(F-1)$.
	This completes the proof.
\end{proof}

Let us remark that it is indeed necessary to separate the case $m=1$ in Theorem~\ref{theo:finite}, since the bound $D \le md(F-1)$ does not hold for $m=1$.
To see this, consider for example the univariate real polynomials $P,Q \in \RR[x]$ given by $P=x^2$ and $Q=x$,
and the set $X=\{-F,\ldots,-1,1,\ldots,F\} \sub \RR$,
so that $|P(X)| = |\{1^2,\ldots,F^2\}|=F$.
Since $0 \notin P(X)$, 
if $Q \equiv_X RP$
then $R(x) = Q(x)/P(x) = 1/x$ for every $x \in X$.
By Lagrange interpolation (originally due to Waring~\cite{Waring79}), the unique univariate polynomial $R(x)$ with $\deg(R) < |X|$ satisfying $R(x)=1/x$ for every $x \in X$ has the maximal degree, $|X|-1$;
explicitly, the coefficient of $x^{|X|-1}$ in $R$ is\footnote{The first equality uses the change of variables $(k,k')=(F+i,F+j)$.} 
\begin{align*}
	\sum_{i \in X} \frac{1}{i}\prod_{j \in X \sm \{i\}}\frac{1}{i-j}
	&= \sum_{\substack{k=0\\k\neq F}}^{2F} \prod_{\substack{k'=0\\k'\neq k}}^{2F} \frac{1}{k-k'}
	= \sum_{\substack{k=0\\k\neq F}}^{2F} \frac{(-1)^k}{(2F)!}\binom{2F}{k}\\
	&= -\frac{(-1)^F}{(2F)!}\binom{2F}{F}
	= \frac{(-1)^{F+1}}{(F!)^2} \neq 0.
\end{align*}
%
Therefore, $\deg(R) = 2F-1$ (which matches~(\ref{eq:bd-m1})) is indeed larger than $2(F-1)$.
The reason that in Theorem~\ref{theo:main} the analogous bound $D \le md(|\F|-1)$ holds even for $m=1$ is that over a finite field $\F$, the univariate function $x \mapsto 1/x$ for $x \neq 0$ 
can actually be interpolated by a polynomial of degree strictly smaller than $|\F|-1$; 
indeed, $x^{|\F|-2} = 1/x$ for every $x \in \F\sm\{0\}$.

\paragraph{Acknowledgments:}
	We thank Shachar Lovett and Robert Robere for helpful discussions.
	We thank Pete L.~Clark for communicating to us the generalization of the main theorem to arbitrary fields (Theorem~\ref{theo:finite}).
	We thank the anonymous referees for useful comments.
	This project was conducted as part of the 2021 New York Discrete Mathematics REU, funded by NSF grant DMS 2051026, and partially supported by award PSC-CUNY TRADB-52-76.

%
%

%
%

\end{document}